\documentclass[12pt]{amsart}

\usepackage[all]{xy}

\newtheorem{theorem}{Theorem}[section]
\newtheorem{corollary}[theorem]{Corollary}
\newtheorem{lemma}[theorem]{Lemma}
\newtheorem{properties}[theorem]{Properties}

\theoremstyle{remark}
\newtheorem{remark}{Remark}

\theoremstyle{definition}
\newtheorem{definition}[theorem]{Definition}
\newtheorem{example}[theorem]{Example}

\begin{document}

\title[Minimal Bratteli diagrams and dimension groups]
{Minimal Bratteli diagrams and dimension groups of AF C$^*$-algebras}
\author{Ryan J. Zerr}
\address{Mathematics Department, University of North Dakota, Grand Forks,
ND, USA} \email{ryan.zerr@und.nodak.edu} 

\begin{abstract}
A method is described which identifies a wide variety of AF
algebra dimension groups with groups of continuous functions.
Since the continuous functions in these groups have domains which
correspond to the set of all infinite paths in what will be called
minimal Bratteli diagrams, it becomes possible, in some cases, to
analyze the dimension group's order preserving automorphisms by
utilizing the topological structure of the associated minimal
Bratteli diagram.
\end{abstract}

\bibliographystyle{amsplain}

\maketitle

\section{Introduction}
\label{section1}

In studying approximately finite-dimensional (AF) C$^*$-algebras
$\displaystyle \mathfrak{A}=\lim_{\longrightarrow}
(\mathfrak{A}_n,\phi_n)$, Bratteli, in his seminal
paper~\cite{bratteli}, introduced a certain infinite graph, now
called a Bratteli diagram, which can be used to encode the nature
of the subalgebras $\mathfrak{A}_n$ and the actions of the
connecting homomorphisms $\phi_n$.  These diagrams, although not
unique for a given algebra, can be used to ascertain certain
characteristics of the algebra, such as in~\cite{bratteli}
and~\cite{lazaretal}.

Subsequently, Elliott~\cite{elliott} proved his theorem showing
that dimension groups provide a complete isomorphism invariant for
AF algebras.  These dimension groups can be realized as a direct
limit of a sequence of scaled ordered groups.  In this context,
the Bratteli diagram is useful since the exact nature of the
groups in this sequence and of the connecting maps between them
can be readily obtained from it.

In the case that the AF algebra under consideration is
commutative, another interesting way in which the Bratteli diagram
plays a role in the description of the dimension group can be
observed. To be specific, let $X$ be a compact metric space with a
basis consisting of sets which are simultaneously open and closed
(clopen).  Then $C(X)$ is AF and the associated dimension group
can be shown to be order isomorphic to the scaled ordered group
$(C(X,\mathbb{Z}),C(X,\mathbb{Z}^+), \chi_{_{X}})$ where
$C(X,\mathbb{Z})$ are the continuous functions $f:X\to\mathbb{Z}$,
$C(X,\mathbb{Z}^+)$ are those functions in $C(X,\mathbb{Z})$ which
only take nonnegative integer values, and $\chi_{_{X}}$ is the
function which is always $1$.

For commutative AF algebras such as this, the spectrum $X$ can be
identified with the set of all infinite paths in the associated
Bratteli diagram.  Therefore, in such a case, the dimension group
can be conveniently described as a set of integer-valued
continuous functions with domain equal to the Bratteli diagram.
This particularly simple description of the dimension group of
$C(X)$ is one source of motivation for the results of this paper.

From the viewpoint of $K$-theory, cf.~\cite{blackadar}, the
dimension group of an AF algebra $\mathfrak{A}$ turns out to be
the $K_0$ group of $\mathfrak{A}$, $K_0(\mathfrak{A})$.  We note
here that ordered groups of the form $C(X,\mathbb{Z})$, where $X$
is a Cantor set, have also made appearances in the $K$-theoretic
calculations of a number of authors.  These have
included~\cite{durandhostskau,giordanoputnamskau,
hermanandputnamandskau,putnam}, where crossed product
C$^*$-algebras were studied, and a determination of their
$K$-theory utilized the fact that $K_0(C(X))$ is isomorphic to
$C(X,\mathbb{Z})$.

In this paper a type of Bratteli diagram, which we refer to as a
\emph{minimal Bratteli diagram} (definition below), will play an
important role. In the case of the AF algebra $C(X)$ above, the
Bratteli diagram is already minimal, and we see that $K_0(C(X))$
can be realized as $C(X,\mathbb{Z})$. More generally, we will
consider certain AF algebras $\mathfrak{A}$ which have proper
minimal diagrams.  That is, the Bratteli diagram itself may not be
minimal, but a subgraph is.  Then, by defining $X_{min}$ to be the
set of all infinite paths in the minimal diagram, it will be
possible to describe $K_0(\mathfrak{A})$ as a subgroup of
$C(X_{min},G)$, where now $G$ is in general a subset of
$\mathbb{Q}$, and under certain hypotheses, $K_0(\mathfrak{A})
\cong C(X_{min},G)$.  In this sense, the results of this paper
provide a direct generalization of the observation that
$K_0(C(X))\cong C(X,\mathbb{Z})$.

For an arbitrary AF algebra, there may correspond multiple,
nonhomeomorphic versions of $X_{min}$.  However, in Section
\ref{section5}, we are able to show that when $\mathfrak{A}$ has a
diagram with a unique graph corresponding to the set $X_{min}$,
then the topological structure of $X_{min}$ provides a way to
discriminate between non-isomorphic algebras (Corollary
\ref{cor}).  This result is a generalization and provides a new
proof of the well known fact that isomorphic commutative AF
algebras have homeomorphic spectra. These facts rely on Theorem
\ref{homeomorphic}, which can also be used to provide information
about the automorphism group of certain dimension groups in terms
of the homeomorphisms on the set $X_{min}$. In particular, those
dimension groups to which this applies will be those arising from
AF algebras that have Bratteli diagrams with a unique choice for
the graph which corresponds to the set $X_{min}$.  It should be
noted that these AF algebras are, in general, different than the
AF algebras with stationary Bratteli diagrams considered
in~\cite{bratteli1, bratteli2, bratteli3}.  There, questions about
stable isomorphisms of these AF algebras and the order preserving
isomorphisms of their corresponding dimension groups are
considered.

In addition to being useful from a theoretical perspective, the
results contained here provide an algorithm for determining the
$K$-theory of many AF algebras, with the resulting dimension
groups being groups of continuous functions.  The types of AF
algebras for which the results of this paper apply include the UHF
and GICAR algebras, with some of the results here being a
generalization of the calculations done for the GICAR algebra
in~\cite{petersandzerr}.  The results of~\cite{petersandzerr} are
in the context of inverse semigroups of partial homeomorphisms,
and it is this dynamical systems perspective on AF algebras which
has, to some extent, motivated the present study. Further recent
work in this area is that of~\cite{donsigandhopenwasser}, which
considers the construction of AF algebras (among others) by
utilizing partial actions.  As pointed out there, there is no
substantive difference between the partial action approach and the
partial homeomorphism approach. As such, the results here might
possibly be generalized to fruitfully study other types of
algebras.

\section{Preliminaries}
\label{section2}

Let $\displaystyle \mathfrak{A}=\lim_{\longrightarrow}
(\mathfrak{A}_n,\phi_n)$ be an AF C$^*$-algebra with
finite-dimensional subalgebras $\mathfrak{A}_n\cong
M_{k(1,n)}\oplus\cdots \oplus M_{k(m_n,n)}$, for all $n\ge 0$
(with $\mathfrak{A}_0\cong \mathbb{C}$), and suppose that for each
$n$, the connecting maps $\phi_n:\mathfrak{A}_n\to
\mathfrak{A}_{n+1}$ are unital injective $*$-homomorphisms. Let
$\mathfrak{D}\subset\mathfrak{A}$ be the abelian subalgebra
$\displaystyle \mathfrak{D}=\lim_{\longrightarrow}
(\mathfrak{D}_n,\phi_n)$, where $\mathfrak{D}_n$ is the subalgebra
of $\mathfrak{A}_n$ spanned by the diagonal matrix units in
$\mathfrak{A}_n$, and let $X_{\mathfrak{A}}$ be the Gelfand
spectrum of $\mathfrak{D}$. Each element of $X_{\mathfrak{A}}$
corresponds uniquely to a decreasing sequence of diagonal matrix
units in $\mathfrak{A}$.  To be more specific, for any $n\ge 0$,
$1\le r\le m_n$, and $1\le s\le k(r,n)$, let
$e_{s,s}(r,n)=e_s(r,n)\in M_{k(r,n)}$ be a diagonal matrix unit of
$\mathfrak{A}$.  Up to unitary equivalence each mapping $\phi_n$
is a standard embedding in the sense of~\cite{power}. Therefore,
$\phi_n(e_s(r,n))$ is a sum of diagonal matrix units in
$\mathfrak{A}_{n+1}$, say
\[
\phi_n(e_s(r,n))=\sum_{l=1}^c e_{\beta_l}(\alpha_l,n+1),
\]
where $c$ depends on the embedding $\phi_n$ and the matrix unit
$e_s(r,n)$.  The elements of $X_{\mathfrak{A}}$ are such that if
$e_{s_0}(r_0,0)\ge e_{s_1}(r_1,1)\ge\ldots$ is one of the
sequences in $X_{\mathfrak{A}}$ then $e_{s_{n+1}}(r_{n+1},n+1)$
appears in the sum $\phi_n(e_{s_n}(r_n,n))$, for all $n\ge 0$.

We now define the set $\widehat{e_s}(r,n)\subset X_{\mathfrak{A}}$
to be the collection of all sequences in $X_{\mathfrak{A}}$ with
the matrix unit $e_s(r,n)$ as the $n$-th coordinate.  Then, the
topology on $X_{\mathfrak{A}}$ has as a basis the collection
\[
\bigvee_{n=1}^\infty \bigvee_{r=1}^{m_n} \bigvee_{s=1}^{k(r,n)}
\widehat{e_s}(r,n).
\]
Of course, from the Gelfand theory, $X_{\mathfrak{A}}$ is compact
Hausdorff, and in fact, each basis element $\widehat{e_s}(r,n)$ is
clopen.  Therefore, $X_{\mathfrak{A}}$ is a $0$-dimensional
compact Hausdorff space.

For the purposes of this paper, a certain closed subset of
$X_{\mathfrak{A}}$ will play an important role. We first define,
for each $n\ge 0$, the set
$X_{min}^n=\bigcup_{r=1}^{m_n}\widehat{e}_1(r,n)$, and let the
subset $X_{min}\subset X_{\mathfrak{A}}$, with the topology it
inherits from $X_{\mathfrak{A}}$, be given by
\[
X_{min}=\bigcap_{n=0}^\infty X_{min}^n.
\]
We will prove the following fact about the set $X_{min}$.

\begin{lemma}
\label{xminlemma}

The sequence $\{X_{min}^n\}_{n=0}^\infty$ is a nested decreasing
sequence of clopen subsets of $X_{\mathfrak{A}}$, and therefore
$X_{min}$ is a nonempty closed subset of $X_{\mathfrak{A}}$.

\end{lemma}

\begin{proof}

For $n\ge 0$ given, let $e_s(r,n)$ be a diagonal matrix unit in
$\mathfrak{A}_n$ where $s>1$ and suppose, to derive a
contradiction, that
\[
\phi_n(e_s(r,n))=\sum_{l=1}^c e_{\beta_l}(\alpha_l ,n+1)
\]
where $\beta_l=1$ for some $l$.  Since (up to unitary equivalence)
$\phi_n$ is a standard embedding, it maps the strictly upper
triangular subalgebra of $\mathfrak{A}_n$ to the strictly upper
triangular subalgebra of $\mathfrak{A}_{n+1}$. Therefore,
$\phi_n(e_{s-1,s}(r,n))$ will be a sum of the form
\[
\phi_n(e_{s-1,s}(r,n))=\sum_{l=1}^c
e_{\gamma_l,\delta_l}(\alpha_l, n+1)
\]
where $\gamma_l< \delta_l$, for all $l$.  The fact that $\phi_n$
is a $*$-homomorphism means
\begin{eqnarray*}
\phi_n(e_s(r,n))&=& \phi_n(e_{s-1,s}
(r,n))^*\phi_n(e_{s-1,s}(r,n))\\
&=& \sum_{l=1}^c
\sum_{l'=1}^c e_{\delta_l,\gamma_l}(\alpha_l,n+1)e_{\gamma_{l'},
\delta_{l'}}(\alpha_{l'},n+1).
\end{eqnarray*}
Because $\phi_n$ maps the diagonal of $\mathfrak{A}_n$ to the
diagonal of $\mathfrak{A}_{n+1}$, if $\gamma_l=\gamma_{l'}$ then
$\delta_l =\delta_{l'}$.  But this implies that $l=l'$.  So
\[
\phi_n(e_s(r,n))=\sum_{l=1}^c e_{\delta_l}(\alpha_l, n+1).
\]
By assumption, it follows that $\delta_{\overline{l}}=1$ for some
$\overline{l}$.  Since $\gamma_l<\delta_l$ for all $l$, we have
that $\gamma_{\overline{l}} < 1$, a contradiction.  Hence, if
$s>1$ then $\phi_n(e_s(r,n))=\sum_{l=1}^c e_{\beta_l}
(\alpha_l,n+1)$ where $\beta_l>1$ for all $l$.

Now, take $(x_1,x_2,\ldots)\in \widehat{e}_1(r,n+1)\subset
X_{min}^{n+1}$.  By the contrapositive of the fact just proven, it
must be that $(x_1,x_2,\ldots)\in \widehat{e}_1(\overline{r},n)$
for some $\overline{r}$.  So, $X_{min}^{n+1} \subset X_{min}^n$,
for all $n\ge 0$.  Hence, $\{X_{min}^n\}_{n=0}^\infty$ is a
nested, decreasing sequence of sets.  Furthermore, since
$\widehat{e}_1(r,n)$ is clopen for all $n$ and all $r$, so too is
$X_{min}^n$ clopen, for all $n\ge 0$.  It then follows that
$X_{min}$ is closed, and by the compactness of $X_{\mathfrak{A}}$
that we have $\bigcap_{n=0}^\infty X_{min}^n \ne\emptyset$.

\end{proof}

At this point we make the observation that for a given AF algebra
$\mathfrak{A}$, the set $X_{min}$ depends on the exact nature of
the standard embeddings $\phi_n$.

\begin{example}
\label{example1}

Consider the AF algebra $\mathfrak{A}$ with Bratteli diagram
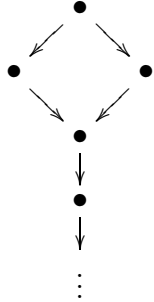
\[
\xymatrix@=12pt{ & \bullet\ar[dl]\ar[dr] & \\
\bullet\ar[dr] & & \bullet\ar[dl] \\
 & \bullet\ar[d] & \\
 & \bullet\ar[d] & \\
 & \vdots &
}\]
Here, we have
\[
X_{\mathfrak{A}}=\{(e_1(1,0),e_1(1,1),\ldots),(e_1(1,0),e_1(2,1),\ldots)
\}.
\]
If we suppose
$\phi_1:\mathfrak{A}_1\cong\mathbb{C}\oplus\mathbb{C} \to
\mathfrak{A}_2\cong M_2$ is such that
\[
\phi_1(a\oplus b)=\left[\begin{array}{cc} a & 0\\ 0 & b
\end{array}\right],
\]
then $X_{min}=\{(e_1(1,0),e_1(1,1),\ldots)\}$.  However, by
assuming
\[
\phi_1(a\oplus b)=\left[\begin{array}{cc} b & 0\\ 0 & a
\end{array}\right],
\]
we arrive at $X_{min}=\{(e_1(1,0),e_1(2,1),\ldots)\}$.

\end{example}

This example is a simple one meant to illustrate the fact that by
changing the sequence of standard embeddings
$\{\phi_n\}_{n=0}^\infty$ to a sequence of unitarily equivalent
standard embeddings, we may arrive at a different set $X_{min}$.
Of course, in this example, these different copies of $X_{min}$
are homeomorphic.  So, in this sense, Example \ref{example1} does
not seem very revealing.  However, we will see in a moment
(Example \ref{nonhomeomorphic}) that by making different,
unitarily equivalent, choices for the sequence of standard
embeddings, non-homeomorphic copies of $X_{min}$ can arise from
the same diagram.

For a given AF algebra $\mathfrak{A}$, the Bratteli diagram
provides a convenient tool for visualizing these different
possibilities  for the set $X_{min}$.  The basis for this begins
with the fact that the set $X_{\mathfrak{A}}$ can be realized as
the set of all infinite paths in the Bratteli diagram for
$\mathfrak{A}$.

To establish some notation for the Bratteli diagram of
$\mathfrak{A}$, we label the vertices at level $n$ by
$k(1,n),\ldots, k(m_n,n)$, analogous to the labelling scheme for
the summands which constitute $\mathfrak{A}_n$. Now, for each
$n\ge 1$, delete enough edges between the vertices
$k(1,n-1),\ldots, k(m_{n-1},n-1)$, and the vertices
$k(1,n),\ldots, k(m_n,n)$, so that each vertex $k(r,n)$, $1\le
r\le m_n$, absorbs exactly one edge from the level above.  The set
of all infinite paths in the resulting subgraph of the Bratteli
diagram will correspond to the set $X_{min}$ for one possible
choice of the standard embeddings $\{\phi_n\}_{n=0}^\infty$.  To
see that this is true, we first note that for all $n\ge 0$, each
vertex $k(i,n)$, $1\le i\le m_n$, can be thought of as
corresponding to an upper left diagonal matrix unit of
$\mathfrak{A}_n$. When viewed from this perspective, the proof of
Lemma \ref{xminlemma} shows us that there is a mapping from the
vertices $\{k(1,n),\ldots, k(m_n,n)\}$ to the vertices
$\{k(1,n-1),\ldots, k(m_{n-1},n-1)\}$, for all $n\ge 1$.
Furthermore, this mapping can be represented by deleting all but
$m_n$ edges between the vertices $\{k(i,n-1):1\le i\le m_{n-1}\}$
and $\{k(j,n):1\le j\le m_n\}$ in the Bratteli diagram so that
each vertex $k(j,n)$ absorbs exactly one edge.  The resulting
subgraph of the Bratteli diagram will then encode, via the set of
all infinite paths, all the possible sequences of the form
$(e_1(i_0,0), e_1(i_1,1),\ldots)$ in $X_{\mathfrak{A}}$, which is
exactly the set $X_{min}$.  By making different choices for the
way we delete the edges, one is in effect making different, albeit
unitarily equivalent, choices for the sequence of standard
embeddings.  In the context of Example \ref{example1}, the two
possible diagrams representing $X_{min}$ would therefore be
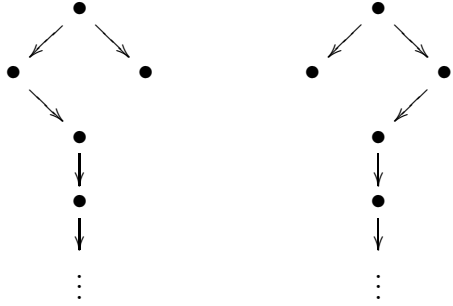
\[
\xymatrix@=12pt{ & \bullet\ar[dl]\ar[dr] & \\
\bullet\ar[dr] & & \bullet \\
 & \bullet\ar[d] & \\
 & \bullet\ar[d] & \\
 & \vdots &
}\hskip 50pt
\xymatrix@=12pt{ & \bullet\ar[dl]\ar[dr] & \\
\bullet & & \bullet\ar[dl] \\
 & \bullet\ar[d] & \\
 & \bullet\ar[d] & \\
 & \vdots &
}
\]

\begin{example}
\label{nonhomeomorphic}

To see that non-homeomorphic copies of $X_{min}$ are possible for
a given Bratteli diagram, consider the GICAR algebra, with
Bratteli diagram
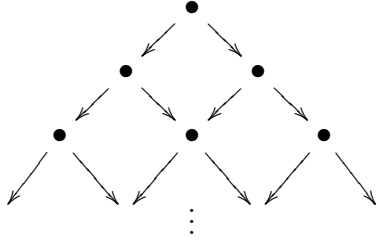
\[
\xymatrix@=12pt{ & & & {\bullet} \ar[dl] \ar[dr] & & &\\
  & & {\bullet} \ar[dl] \ar[dr] & & {\bullet} \ar[dl] \ar[dr] & &\\
  & {\bullet} \ar[dl] \ar[dr] & & {\bullet} \ar[dl] \ar[dr] & &
  {\bullet} \ar[dl] \ar[dr] &\\
  & & & {\vdots} & & & }
\]
Two possible diagrams representing $X_{min}$, achieved through
edge deletions as described above, are
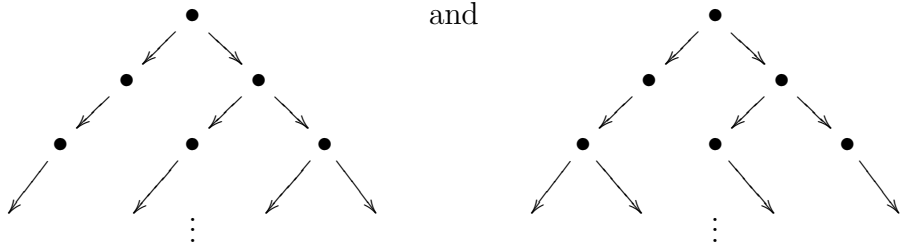
\[
\xymatrix@=12pt{ & & & {\bullet} \ar[dl] \ar[dr] & & &\\
  & & {\bullet} \ar[dl] & & {\bullet} \ar[dl] \ar[dr] & &\\
  & {\bullet} \ar[dl]  & & {\bullet} \ar[dl]  & &
  {\bullet} \ar[dl] \ar[dr] &\\
  & & & {\vdots} & & & }\hskip 10pt \mbox{ and } \hskip 10pt
  \xymatrix@=12pt{ & & & {\bullet} \ar[dl] \ar[dr] & & &\\
  & & {\bullet} \ar[dl] & & {\bullet} \ar[dl] \ar[dr] & &\\
  & {\bullet} \ar[dl]\ar[dr]  & & {\bullet} \ar[dr]  & &
  {\bullet} \ar[dr] &\\
  & & & {\vdots} & & & }
\]
The sets of all infinite paths in these diagrams are homeomorphic
to
\[
\overline{\left\{\frac{1}{n}:n\in\mathbb{Z}^+\right\}}\hskip 25pt
\mbox{ and } \hskip 25pt \overline{\left\{\pm\left(
1-\frac{1}{n}\right):n\in \mathbb{Z}^+\right\}},
\]
respectively.  Since the former has one non-isolated point and the
latter two such points, they are clearly not homeomorphic.

\end{example}

\section{Minimal Bratteli Diagrams}

Given an AF algebra $\displaystyle
\mathfrak{A}=\lim_{\longrightarrow} (\mathfrak{A}_n,\phi_n)$,
with, for all $n\ge 0$, $\mathfrak{A}_n\cong
M_{k(1,n)}\oplus\cdots\oplus M_{k(m_n,n)}$, consider the sequence
$\{ m_n\}_{n=0}^\infty$ of positive integers.  Two possibilities
exist:
\begin{itemize}
\item[(1)] $\displaystyle \limsup_{n\to\infty} m_n=L<\infty$; or
\item[(2)] $\displaystyle \limsup_{n\to\infty} m_n=\infty$.
\end{itemize}
For case (1), there exists a subsequence
$\{m_{n_k}\}_{k=1}^\infty$ of $\{m_n\}_{n=0}^\infty$ such that
$m_{n_k}=L$, for all $k\ge 1$, and for case (2), there exists a
subsequence $\{m_{n_k}\}_{k=1}^\infty$ such that
$m_{n_k}<m_{n_{k+1}}$, for all $k\ge 1$.  For either situation,
the sequence $\{\mathfrak{A}_n\}_{n=0}^\infty$ of
finite-dimensional subalgebras can be \emph{contracted} to the
subsequence $\{\mathfrak{A}_{n_k}\}_{k=0}^\infty$ (where $n_0=0$)
with
\[
\mathfrak{A}=\lim_{\longrightarrow} (\mathfrak{A}_n,\phi_n)
=\lim_{\longrightarrow} (\mathfrak{A}_{n_k},\phi_{n_k}).
\]
Therefore, without a loss of generality, we may assume that every
AF algebra $\mathfrak{A}$ is such that either
$\{m_n\}_{n=0}^\infty$ is a constant sequence or it is strictly
monotonically increasing.

\begin{remark}

There may be many different ways to contract (\emph{telescope} in
the terminology of~\cite{giordanoputnamskau}) a given sequence
$\{\mathfrak{A}_n\}_{n=0}^\infty$ such that the above is true of
$\{m_n\}_{n=0}^\infty$.  In fact, it may be possible to find
multiple contractions of both types.

\end{remark}

As described by~\cite{davidson} and~\cite{power}, each connecting
homomorphism $\phi_n:\mathfrak{A}_n\to\mathfrak{A}_{n+1}$, for an
AF algebra $\displaystyle \mathfrak{A}=\lim_{\longrightarrow}
(\mathfrak{A}_n,\phi_n)$, has associated with it a
{\emph{multiplicity matrix}} $\overline{A}_{n,n+1}\in
M_{m_{n+1},m_n}(\mathbb{N})$ which, up to unitary equivalence,
completely describes the action of the mapping.  Assuming that the
sequence $\{\mathfrak{A}_n\}_{n=0}^\infty$ has been contracted so
as to fit into one of cases (1) or (2), we will, for the purposes
of this paper, \emph{assume that each of the multiplicity matrices
$\overline{A}_{n,n+1}$ has full rank}.  That is, for all $n\ge 0$,
we will assume $\mbox{rank}(\overline{A}_{n,n+1})=m_n$.  So, we
are interested in those AF algebras for which there exists a
contraction of either type (1) or (2) and for which this
contraction results in all multiplicity matrices being of full
rank.

For the moment, we focus our attention on case (2), where $m_n<
m_{n+1}$, for all $n\ge 0$.  We want to show that this case can be
reduced to the situation where $m_n=n+1$, for all $n\ge 0$.

Suppose that for some $n\ge 0$ we have $m_n+1<m_{n+1}$.  Since we
are assuming $\mbox{rank}(\overline{A}_{n,n+1})=m_n$, there are
$m_n$ linearly independent rows in $\overline{A}_{n,n+1}$.  Let
$P\in M_{m_{n+1}}$ be the permutation matrix which permutes the
rows of $\overline{A}_{n,n+1}$ in such a way that
$P\overline{A}_{n,n+1} =[b_{ij}]\in M_{m_{n+1},m_n}(\mathbb{N})$
has $m_n$ linearly independent rows of $\overline{A}_{n,n+1}$ as
its last $m_n$ rows. Then, the matrices
\[
B_1=\left[\begin{array}{ccc} b_{1,1} & \cdots & b_{1,m_n}\\
1 & & \\
  & \ddots & \\
  & & 1
\end{array}\right],
\]
\[
B_2=\left[\begin{array}{cccc} 1 & 0 & \cdots & 0\\
0 & b_{2,1} & \cdots & b_{2,m_n}\\
0 & 1 & & \\
\vdots & & \ddots & \\
0 & & & 1
\end{array}\right],\ldots,
\]
\[
B_{m_{n+1}-m_n-1}=\left[\begin{array}{cccccc} 1 & & & 0 & \cdots &
0\\
 & \ddots & & \vdots & & \vdots\\
 & & 1 & 0 & \cdots & 0\\
0&\cdots &0& b_{m_{n+1}-m_n-1,1}& \cdots & b_{m_{n+1}-m_n-1,m_n}\\
0&\cdots &0& 1 &  & \\
\vdots& &\vdots&   & \ddots & \\
0&\cdots &0&   & & 1
\end{array}\right],
\]
and
\[
B_{m_{n+1}-m_n}=\left[\begin{array}{cccccc} 1 & & & 0 & \cdots &
0\\
 & \ddots & & \vdots & & \vdots\\
 & & 1 & 0 & \cdots & 0\\
0&\cdots &0& b_{m_{n+1}-m_n,1}& \cdots & b_{m_{n+1}-m_n,m_n}\\
0&\cdots &0& b_{m_{n+1}-m_n+1,1} &\cdots & b_{m_{n+1}-m_n+1,m_n} \\
\vdots& &\vdots&\vdots   & &\vdots \\
0&\cdots &0& b_{m_{n+1},1} &\cdots & b_{m_{n+1},m_n}
\end{array}\right]
\]
are such that
\[
P\overline{A}_{n,n+1}=B_{m_{n+1}-m_n}B_{m_{n+1}-m_n-1}\cdots B_2
B_1,
\]
and each matrix $B_i$, $1\le i\le m_{n+1}-m_n$, has full rank.
Furthermore, each matrix $B_1,\ldots, B_{m_{n+1}-m_n-1}, P^{-1}
B_{m_{n+1}-m_n}$ can be viewed as a multiplicity matrix for some
unital embedding between finite-dimensional C$^*$-algebras.  We
can therefore \emph{dilate} (\emph{microscope} in the terminology
of~\cite{giordanoputnamskau}) the sequence
$\{\mathfrak{A}_n\}_{n=0}^\infty$ to a sequence
$\{\mathfrak{B}_n\}_{n=0}^\infty$ of finite-dimensional
C$^*$-algebras such that $\{\mathfrak{B}_n\}_{n=0}^\infty$ can be
contracted to a subsequence equal to
$\{\mathfrak{A}_n\}_{n=0}^\infty$.  Hence, we would have
\[
\mathfrak{A}=\lim_{\longrightarrow} (\mathfrak{A}_n,\phi_n)
=\lim_{\longrightarrow} (\mathfrak{B}_n,\phi'_n),
\]
with the additional properties that the sequence
$\{\mathfrak{B}_n\}_{n=0}^\infty$ is such that $m_n=n+1$, for all
$n\ge 0$, and the multiplicity matrices corresponding to the
connecting embeddings all have full rank. We therefore assume,
without a loss of generality, that for case (2), our AF algebras
$\mathfrak{A}$ are such that $m_n=n+1$, for all $n\ge 0$.

As mentioned in Section \ref{section2}, for a given AF algebra,
many different possibilities may exist for the set $X_{min}$.  For
the proof of this paper's main result, it will be necessary for
$X_{min}$ to be chosen so that it is, in some sense, large.  Our
assumption that each multiplicity matrix has full rank is
sufficient to guarantee such a choice always exists, and the
remainder of this section will be devoted to a consideration of
precisely what we mean by the word ``large''.

To begin, we consider an AF algebra of type (2).  That is,
$\displaystyle \mathfrak{A}=\lim_{\longrightarrow}
(\mathfrak{A}_n,\phi_n)$ where $m_n=n+1$ and
$\mbox{rank}(\overline{A}_{n,n+1})=n+1$, for all $n\ge 0$.  We
will employ the method described in Section \ref{section2} of
deleting certain edges between the vertices in each level of the
Bratteli diagram of $\mathfrak{A}$ in order to describe $X_{min}$
as a set of infinite paths in the resulting subgraph.

For the moment, we adopt a general perspective.  Let $n\ge 1$ be
given.  We will consider graphs with the following properties:
\begin{properties}
\label{properties}
\end{properties}
\begin{itemize}

\item[(I)] There are $2n+1$ vertices which are arranged so that
$n$ vertices appear in one horizontal row (which we will refer to
as level 1) and the remaining $n+1$ vertices appear in a
horizontal row below the first (which we will refer to as level
2).

\item[(II)] The only edges are those connecting vertices at
different levels.  In other words, no vertices at the same level
are connected by an edge.

\item[(III)] Every vertex is connected to another by at least on
edge.

\end{itemize}
We label the set of all such graphs $G_n$.

\begin{remark}

The AF algebras with $m_n=n+1$, for all $n\ge 0$, are built from
graphs of this form.

\end{remark}

\begin{definition}

Given a graph $\Gamma\in G_n$, we will call the graph $\gamma$ a
\emph{reduction of} $\Gamma$ if $\gamma$ is a subgraph of $\Gamma$
obtained by deleting only edges and $\gamma\in G_n$.  A graph in
$G_n$ will be called \emph{minimal} if it has $n+1$ edges.

\end{definition}

\begin{remark}

Any graph in $G_n$ must have at least $n+1$ edges by Properties
\ref{properties} (III).  Therefore, no nontrivial reductions of
minimal graphs exist.

\end{remark}

\begin{example}

It is easy to see that minimal reductions are not unique.  The
graph
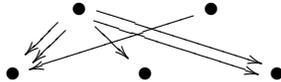
\[
\xymatrix@=12pt{ & {\bullet}\ar@<2pt>[dl]\ar@<-2pt>[dl]\ar[dr]
\ar@<2pt>[drrr]\ar@<-2pt>[drrr] & & {\bullet}\ar[dlll] & \\
 {\bullet} & & {\bullet} & & {\bullet}
  }
\]
is an element of $G_2$, with both
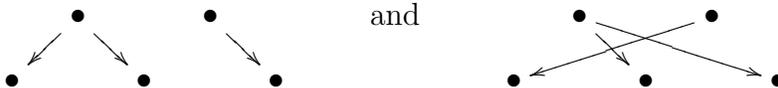
\[
\xymatrix@=12pt{
 & {\bullet}\ar[dl]\ar[dr] & & {\bullet}\ar[dr] &\\
 {\bullet} & & {\bullet} & & {\bullet}
}\hskip 25pt \mbox{ and }\hskip 25pt \xymatrix@=12pt{
 & {\bullet}\ar[dr]\ar[drrr] & & {\bullet}\ar[dlll] &\\
 {\bullet} & & {\bullet} & & {\bullet}
}
\]
being minimal reductions.

\end{example}

\begin{example}
\label{nonminreduct}

For an arbitrary $n$, not every element of $G_n$ has a minimal
reduction. The following element of $G_3$ is such that every way
of deleting all but four edges results in an element not in $G_3$.
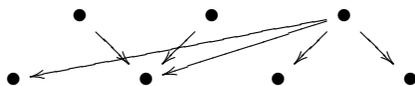
\[
\xymatrix@=12pt{
 & {\bullet}\ar[dr] & & {\bullet}\ar[dl] & & {\bullet}
 \ar[dlllll]\ar[dlll]\ar[dl]\ar[dr] &\\
{\bullet} & & {\bullet} & & {\bullet} & & {\bullet} }
\]

\end{example}

For our purposes, we will be interested in Bratteli diagrams
which, at each level, have minimal reductions. Example
\ref{nonminreduct} demonstrates that not all Bratteli diagrams
will have minimal reductions at each level.  However, if a
Bratteli diagram does have a minimal reduction at each level, then
by carrying out this reduction for all levels in the diagram, the
corresponding ``reduction'' of the Bratteli diagram will
correspond to one possible choice for $X_{min}$ and will itself be
a Bratteli diagram.  One might call such a graph a \emph{minimal
Bratteli diagram} since deleting any more edges will create a
subgraph which is no longer a Bratteli diagram. This will be an
important object for our main result. We are therefore interested
in answering the following somewhat more general question:
\begin{itemize}

\item For any $n\ge 1$, how can we decide which elements of $G_n$
have a minimal reduction (in $G_n$)?

\end{itemize}

At this point we will begin to make use of the assumption we have
made about the multiplicity matrices $\overline{A}_{n,n+1}$.
Namely, that each has full rank.  Under such circumstances we can
guarantee that a graph $\Gamma\in G_n$ has a minimal reduction for
all $n\ge 1$.

\begin{theorem}
\label{graph}

Let $\Gamma\in G_n$ for any $n\ge 1$ and suppose $M_\Gamma$ is the
multiplicity matrix which describes $\Gamma$.  If
$\mbox{\textup{rank}}(M_\Gamma)=n$ then $\Gamma$ has a minimal
reduction.

\end{theorem}

\begin{proof}

We will proceed via induction on $n$.  The induction basis is
provided by the case $n=1$, where it is easy to see that all
graphs $\Gamma\in G_1$ have minimal reductions.

Next, suppose the result holds for $n$ and let $\Gamma\in
G_{n+1}$. We will assume $\mbox{rank}(M_\Gamma)=n+1$ and write
$M_\Gamma= [a_{ij}]$, $1\le i\le n+2$, $1\le j\le n+1$.  Before
proceeding further we remark that by permuting the rows and
columns of $M_\Gamma$ we are merely rearranging the vertices in
the graph. For example, a row permutation amounts to rearranging
the vertices at level $2$ and a column permutation amounts to
rearranging the vertices at level $1$.  Of course, if it is
possible to obtain a minimal reduction of this permuted form of
the original graph, then by reversing the permutations, we will
also have a minimal reduction of the original graph.  We will
therefore work with various matrices obtained from $M_\Gamma$
through row and column permutations in order to obtain our result.

We consider two cases:
\begin{itemize}

\item[(a)] There exists $j_0\in\{ 1,\ldots,n+1\}$ such that the
submatrix $[a_{ij}]$, $1\le i\le n+2$, $1\le j\le n+1$, $j\ne
j_0$, has only nonzero rows.

\item[(b)] For every $j_0\in\{ 1,\ldots,n+1\}$, the submatrix
$[a_{ij}]$, $1\le i\le n+2$, $1\le j\le n+1$, $j\ne j_0$, has a
zero row.

\end{itemize}

Consider case (a).  Permute the columns of $M_\Gamma$ so that the
submatrix that results from omitting the last column of the
permuted matrix has only nonzero rows.  For notational
convenience, we will continue to write $M_\Gamma$ and $[a_{ij}]$
for these permuted forms of the original multiplicity matrix. Now,
permute the rows so that the top $n+1$ rows are linearly
independent, which we can do since $\mbox{rank}(M_\Gamma)=n+1$.
Thus, the submatrix $[a_{ij}]$, $1\le i,j\le n+1$, is nonsingular.

Because $[a_{ij}]$, $1\le i,j\le n+1$, is nonsingular, the matrix
$[a_{ij}]$, $1\le i\le n+1$, $1\le j\le n$, has rank $n$.  Thus,
there exists $i$, which, without a loss of generality, we may
suppose equals $n+1$, such that the submatrix $[a_{ij}]$, $1\le
i,j\le n$, is nonsingular.

At this point we make an assumption that is justified later (Lemma
\ref{rowdel}). Assume $a_{n+1,n+1}\ne 0$.  Permute the last two
rows of $M_\Gamma$ to finally arrive at a matrix with the
following characteristics:
\begin{itemize}

\item[(i)] Omitting the $(n+1)$-st column leaves all nonzero rows;

\item[(ii)] The first $n$ rows of the submatrix $[a_{ij}]$, $1\le
i\le n+1$, $1\le j\le n$, are linearly independent; and

\item[(iii)] The entry $a_{n+2,n+1}$ is nonzero.

\end{itemize}
Since (i) and (ii) hold, the submatrix $[a_{ij}]$, $1\le i\le
n+1$, $1\le j\le n$, is the multiplicity matrix for a graph in
$G_n$ with rank $n$.  Thus, by the induction hypothesis, there
exists a minimal reduction of this graph (which is just a subgraph
of $\Gamma$).  Since (iii) holds, there is at least one edge
connecting the last vertex at level $1$ with the last vertex at
level $2$.  By deleting all other edges which connect these last
two vertices to any others, we obtain a minimal reduction of
$\Gamma$.

Next, consider case (b).  The assumption implies that at least
$n+1$ rows have exactly one nonzero entry and for different rows,
the columns in which these entries appear are different. We may
assume without a loss of generality that these rows are the first
$n+1$ rows. So, the first $n+1$ vertices at level $2$ are
connected to exactly one vertex at level $1$, and different level
$2$ vertices are connected to different level $1$ vertices.
Finally, since the last row is not zero, the last vertex at level
$2$ is connected to at least one vertex at level $1$. Thus, a
minimal reduction in $G_{n+1}$ is possible.

\end{proof}

We now justify an assumption made in the proof of the previous
theorem.

\begin{lemma}
\label{rowdel}

Given an invertible matrix $B=[b_{ij}]\in M_n(\mathbb{C})$, there
exists $k$, $1\le k\le n$, such that the submatrix $[b_{ij}]$,
$1\le i,j\le n$, $i\ne k$, $j\ne n$, is nonsingular and
$b_{k,n}\ne 0$.

\end{lemma}

\begin{proof}

First, if there exists $1\le k\le n$ such that $b_{k,1}=\cdots
=b_{k,n-1}=0$, then it must be that $b_{k,n}\ne 0$. Furthermore,
since the dimension of the set
$\mbox{span}\{[b_{i,1},\ldots,b_{i,n-1}]: 1\le i\le n\}$ is $n-1$,
the desired result is achieved.  To complete the proof we consider
the case where $[b_{i,1},\ldots, b_{i,n-1}]\ne 0$, for every $1\le
i\le n$.

Assume, without a loss of generality, that the first $i_0$ rows of
the matrix $B$ end in $0$ (i.e., $b_{i,n}=0$, for all $1\le i\le
i_0$) and the remaining rows end in a nonzero number (i.e.,
$b_{i,n}\ne 0$, for all $i_0+1\le i\le n$).  Of course, $i_0<n$
since $B$ is nonsingular.  If we assume that there exist scalars
$\alpha_1,\ldots, \alpha_{i_0}$ such that
\[
\alpha_1[b_{1,1},\ldots,b_{1,n-1}]+\cdots+ \alpha_{i_0}
[b_{i_0,1},\ldots,b_{i_0,n-1}]=0
\]
then
\[
\alpha_1[b_{1,1},\ldots,b_{1,n-1},b_{1,n}]+\cdots+ \alpha_{i_0}
[b_{i_0,1},\ldots,b_{i_0,n-1},b_{i_0,n}]=0
\]
as well since $b_{1,n}=\cdots=b_{i_0,n}=0$.  Because these later
vectors are linearly independent, it must be that $\alpha_1=\cdots
=\alpha_{i_0}=0$, and therefore, the set
\[
\{[b_{i,1},\ldots, b_{i,n-1}]:1\le i\le i_0\}
\]
is linearly independent.

We know the set $\{[b_{i,1},\ldots,b_{i,n-1}]: 1\le i\le n\}$ is
linearly dependent, and so there exist coefficients
$\alpha_1,\ldots,\alpha_n$, not all zero, such that
\[
\sum_{i=1}^n\alpha_i[b_{i,1},\ldots,b_{i,n-1}]=0.
\]
Furthermore, since $\{[b_{i,1},\ldots,b_{i,n-1}]: 1\le i\le i_0\}$
is a linearly independent set, it must be that at least one
coefficient $\alpha_{i_1}$, such that $i_1>i_0$, is nonzero.  But
then,
\begin{eqnarray*}
\mathbb{C}^{n-1}&=& \mbox{span}\{ [b_{i,1},\ldots,b_{i,n-1}]: 1\le
i\le n\}\\
&=& \mbox{span}\{[b_{i,1},\ldots,b_{i,n-1}]:1\le i\le n, i\ne
i_1\},
\end{eqnarray*}
which implies $\{[b_{i,1},\ldots,b_{i,n-1}]:1\le i\le n,i\ne
i_1\}$ is linearly independent.  Hence, the submatrix $[b_{ij}]$,
$1\le i\le n$, $i\ne i_1$, $1\le j\le n-1$, is nonsingular and
$b_{i_1,n}\ne 0$.

\end{proof}

This proves that for those AF algebras which are of type (2), a
minimal Bratteli diagram corresponding to the set $X_{min}$ exists
under the assumption that each multiplicity matrix has full rank.
For algebras of type (1), we have an analogous result.

\begin{theorem}
\label{finitegraphs}

Suppose $\Gamma$ is a graph consisting of $n$ vertices at both
levels $1$ and $2$ and \textup{(II)} and \textup{(III)} of
Properties \textup{\ref{properties}} are satisfied. If the
multiplicity matrix $M_\Gamma$ which describes the edges has rank
$n$, then $\Gamma$ has a minimal reduction \textup{(}a reduction
with $n$ edges\textup{)} which satisfies \textup{(II)} and
\textup{(III)} of Properties \textup{\ref{properties}}.

\end{theorem}

\begin{proof}

The proof is very similar to that of Theorem \ref{graph}.

\end{proof}

The results obtained here will be used in the subsequent section
to prove a result about the $K_0$ groups of a large class of AF
algebras.  Despite that motive for their inclusion, they are
interesting in their own right. After all, being able to reduce a
Bratteli diagram in the way described here means that a
sub-Bratteli diagram exists which is minimal in some sense.
Specifically, deleting any more edges will result in a subgraph
which is no longer itself a Bratteli diagram. Of course, all
Bratteli diagrams are reducible to sub-Bratteli diagrams (possibly
in a trivial way). However, the reductions here to the level of
$X_{min}$ are as far as one can go.  Deleting any more edges,
without deleting any vertices, will result in a subgraph which is
no longer a Bratteli diagram.

\section{Dimension Groups and  Minimal Bratteli Diagrams}
\label{sec4}

As is well known, if $X$ is a $0$-dimensional (basis consisting of
clopen sets) compact metric space, then there exists a sequence
$\{ E_n\}_{n=0}^\infty$ of successively finer partitions of $X$
which generate the topology and consist of clopen sets. Therefore,
$\{ C(E_n)\}_{n=0}^\infty$, where $C(E_n)$ consists of those
functions constant on elements of the partition $E_n$, is an
increasing sequence of finite-dimensional C$^*$-algebras
(isomorphic to $\mathbb{C}^{|E_n|}$).  Since
$C(X)=\overline{\bigcup_{n\ge 0} C(E_n)}$, it follows that $C(X)$
is AF.

For any $n\ge 0$, the dimension group $K_0(C(E_n))$ of $C(E_n)$ is
easily seen to be isomorphic to $\left( C(E_n,\mathbb{Z}),
C(E_n,\mathbb{Z}^+), \chi_{_{X}}\right)$, where
$C(E_n,\mathbb{Z})$ are the continuous functions from $X$ to
$\mathbb{Z}$ constant on the elements of $E_n$.  As such, one can
conclude that
\[
K_0(C(X))=\lim_{\longrightarrow} K_0(C(E_n))\cong \left(
C(X,\mathbb{Z}), C(X,\mathbb{Z}^+),\chi_{_{X}}\right).
\]

One aspect of this well known example which we want to draw
attention to is that $X_{min}=X$. Therefore, in a natural way, the
dimension group of $C(X)$ can be realized as a group of continuous
functions on $X_{min}$. In the context of AF
groupoids,~\cite{petersandzerr} presents additional examples for
which this also holds true. We intend in this section to
generalize the results of~\cite{petersandzerr} in order to
demonstrate that the dimension groups of those AF algebras in a
certain class can be realized in this same way, as groups of
continuous functions on $X_{min}$.

Here, as in the previous section, we begin with an AF algebra
$\mathfrak{A}$ such that the sequence $\{m_n\}_{n=0}^\infty$ is
strictly monotonically increasing and each multiplicity matrix
$\overline{A}_{n,n+1}$ has full rank.  Of course, by one of our
results, we may assume without a loss of generality that
$m_n=n+1$, for all $n\ge 0$.

By Theorem \ref{graph} we can reduce the Bratteli diagram
representing $\mathfrak{A}$ to a sub-Bratteli diagram whose set of
all infinite paths corresponds to $X_{min}$, as in the discussion
of Section \ref{section2}.  Put in somewhat more concrete terms,
Theorem \ref{graph} tells us that the standard embeddings can be
chosen so that for each $n\ge 1$, there exist integers $1\le r'_n<
r_n\le n+1$ and $1\le a_n\le n$ with
\[
\widehat{e}_1(r'_n,n)\cup\widehat{e}_1(r_n,n)\subset
\widehat{e}_1(a_n,n-1),
\]
and that there exists a bijection
\[
\sigma:\{i:1\le i\le n+1,i\ne r'_n,r_n\} \to \{j:1\le j\le n, j\ne
a_n\}
\]
such that $\widehat{e}_1(i,n)\subset\widehat{e}_1(\sigma(i),n-1)$.

To achieve our result, we will now define, for all $n\ge 1$, a
linear map $R_n:\mathbb{C}^{n+1}\to C(X_{min},\mathbb{C})$ by
\[
R_n(\alpha_1,\ldots,\alpha_{n+1})=\sum_{l=0}^n \alpha_{l+1}
\chi_{_{B(r_l,l)}},
\]
where $r_0=1$ and, as a notational convenience, we write $B(i,n)$
to denote the set $\widehat{e}_1(i,n)\cap X_{min}$.  Note that we
are now assuming that a specific choice for $X_{min}$ has been
made.  Thus, the sequences $\{r_n\}_{n=1}^\infty$ and
$\{r'_n\}_{n=1}^\infty$ are fixed.  The following technical result
about the maps $R_n$ will be useful to us in what follows.

\begin{lemma}
\label{functionform}

For all $n\ge 1$, if $(\alpha_1,\ldots,\alpha_{n+1})^T\in
\mathbb{C}^{n+1}$ then there exists $(\beta_1,\ldots,
\beta_{n+1})^T\in\mathbb{C}^{n+1}$ such that
\[
R_n(\beta_1,\ldots,\beta_{n+1})=\sum_{l=1}^{n+1} \alpha_l
\chi_{_{B(l,n)}}.
\]

\end{lemma}

\begin{proof}

For $n=1$, set $f=\alpha_1\chi_{_{B(1,1)}}+ \alpha_2
\chi_{_{B(2,1)}}$.  Let $\beta_1=\alpha_1$ and
$\beta_2=\alpha_2-\alpha_1$. Then,
\begin{eqnarray*}
R_1(\beta_1,\beta_2)&=& \beta_1\chi_{_{B(1,0)}}
+\beta_2\chi_{_{B(2,1)}}\\
&=& \alpha_1\chi_{_{B(1,0)}}+ (\alpha_2-\alpha_1)
\chi_{_{B(2,1)}}.
\end{eqnarray*}
We note that $B(1,1)\cup B(2,1)= B(1,0)$, and so,
\[
R_1(\beta_1,\beta_2)= \alpha_1\chi_{_{B(1,1)}}
+\alpha_1\chi_{_{B(2,1)}} +(\alpha_2-\alpha_1) \chi_{_{B(2,1)}}=f.
\]

Now, assume that for every $f=\sum_{l=1}^{n+1} \alpha_l
\chi_{_{B(l,n)}}$ there exists a vector $(\beta_1,\ldots,
\beta_{n+1})^T$ such that $R_n(\beta_1,\ldots, \beta_{n+1}) =f$,
where the $\beta_i$ are linear combinations of the $\alpha_j$. Let
$g=\sum_{l=1}^{n+2} \gamma_l\chi_{_{B(l,n+1)}}$.  By our earlier
comments, there exist $1\le r'_{n+1}<r_{n+1}\le n+2$ and $1\le
a_{n+1}\le n+1$ such that $B(r'_{n+1},n+1)\cup B(r_{n+1},n+1)=
B(a_{n+1},n)$. Furthermore, there exists a bijection $\sigma$ such
that $B(i,n+1)= B(\sigma(i),n)$. Then,
\begin{eqnarray*}
g&=& \sum_{l=1}^{n+2}\gamma_l\chi_{_{B(l,n+1)}} = \sum_{l=1, l\ne
r_{n+1},r'_{n+1}}^{n+2} \gamma_l
\chi_{_{B(\sigma(l),n)}}\\
&+&\gamma_{r_{n+1}} \chi_{_{B(r_{n+1},n+1)}} +\gamma_{r'_{n+1}}
\chi_{_{B(r'_{n+1},n+1)}}\\
&=& \sum_{l=1,l\ne r_{n+1}}^{n+2} \gamma_l
\chi_{_{B(\sigma(l),n)}} +(\gamma_{r_{n+1}}-\gamma_{r'_{n+1}})
\chi_{_{B(r_{n+1},n+1)}}
\end{eqnarray*}
where we define $\sigma(r'_{n+1})=a_{n+1}$.  By the induction
hypothesis, there exists $(\beta_1,\ldots,\beta_{n+1})$ such that
\[
\sum_{l=1,l\ne
r_{n+1}}^{n+2}\gamma_l\chi_{_{B(\sigma(l),n)}}=R_n(\beta_1,\ldots,\beta_{n+1}).
\]
Thus,
\begin{eqnarray*}
g&=&
R_n(\beta_1,\ldots,\beta_{n+1})+(\gamma_{r_{n+1}}-\gamma_{r'_{n+1}})
\chi_{_{B(r_{n+1},n+1)}}\\
&=& \sum_{l=0}^n \beta_{l+1}\chi_{_{B(r_l,l)}}
+(\gamma_{r_{n+1}}-\gamma_{r'_{n+1}})\chi_{_{B(r_{n+1},n+1)}}\\
&=& \sum_{l=0}^{n+1}\beta_{l+1}\chi_{_{B(r_l,l)}} =R_{n+1}
(\beta_1,\ldots, \beta_{n+2}),
\end{eqnarray*}
where $\beta_{n+2}=\gamma_{r_{n+1}}-\gamma_{r'_{n+1}}$.  Hence, by
induction, the result holds for all $n\ge 1$.

\end{proof}

We note that $R_n$ injects for all $n\ge 1$. This follows from the
fact that by construction, the vector $(\beta_1,\ldots,
\beta_{n+1})^T$ which accomplishes $R_n(\beta_1, \ldots,
\beta_{n+1})=0$ has components which are linear combinations of
zeros.  Therefore, $\beta_1=\cdots= \beta_{n+1} =0$.

Now since, by assumption, $\overline{A}_{n,n+1}$ has full rank,
for all $n\ge 0$, it is possible, by adding appropriate columns
$[a_{1,n+2},\ldots,a_{n+2,n+2}]^T\in\mathbb{N}^{n+2}$ to create a
sequence $\{A_{n,n+1}\}_{n=0}^\infty$ of nonsingular matrices
where
\[
A_{n,n+1}=\left[\begin{array}{cc} & a_{1,n+2}\\
\overline{A}_{n,n+1} & \vdots\\
 & a_{n+2,n+2}\end{array}\right]\in M_{n+2}(\mathbb{N}).
\]

\begin{remark}

At this point, any choice for $[a_{1,n+2},\ldots,a_{n+2,n+2}]^T$
which makes $A_{n,n+1}$ nonsingular is appropriate.  In fact,
there is no \emph{a priori} reason that the $a_{i,n+2}$ can not be
elements of $\mathbb{C}$. However, as we will see in a moment, in
certain instances a somewhat more restrictive choice will be
desirable.

\end{remark}

To compute the dimension group $K_0(\mathfrak{A})$, we utilize the
fact that $\displaystyle K_0(\mathfrak{A})=\lim_{\longrightarrow}
K_0(\mathfrak{A}_n)$.  However, to be more explicit about the
nature of $\displaystyle \lim_{\longrightarrow}
K_0(\mathfrak{A}_n)$, we first define, for all $n\ge 1$,
\[
A_n=\left[A_{0,1}^{-1}\oplus I_{n-1}\right]\cdots
\left[A_{n-2,n-1}^{-1}\oplus I_1\right]A_{n-1,n}^{-1}.
\]
Then, let $\Phi_n:\mathbb{Z}^{n+1}\to C(X_{min},G)$ be given by
$\Phi_n=R_n\circ A_n$, for all $n\ge 1$, where $G$ is an abelian
group that will now be described.  Of course, due to the
definition of $R_n$, the nature of $G$ depends on the range of
$A_n$ when its domain is taken as $\mathbb{Z}^{n+1}$.

Letting $n\ge 1$ be given, it is clear that
\[
A_n^{-1}=A_{n-1,n}\left[A_{n-2,n-1}\oplus I_1\right] \cdots
\left[A_{1,2}\oplus I_{n-2}\right]\left[ A_{0,1}\oplus
I_{n-1}\right]
\]
and, by the multiplicativity of the determinant, that
\[
\left|A_n^{-1}\right| =\left| A_{n-1,n}\right|\cdot \left|
A_{n-2,n-1}\right|\cdot\cdots\cdot\left|A_{1,2}\right|\cdot
\left|A_{0,1}\right|.
\]
As, for example, in~\cite[pages 20--21]{hornandjohnson}, we can
then write
\[
A_n=\frac{1}{\left|A_n^{-1}\right|}\mbox{adj}
\left(A_n^{-1}\right),
\]
where $\mbox{adj}(A_n^{-1})$ is the \emph{adjugate} or
\emph{classical adjoint} of the matrix $A_n^{-1}$.  In defining
each matrix $A_{n,n+1}$, if we choose columns of integers, then
each matrix $A_n^{-1}$ will be a product of matrices with integer
entries, and therefore itself must have integer entries. But, by
the definition of $\mbox{adj}(A_n^{-1})$, it will then follow that
$\mbox{adj}(A_n^{-1})$ has integer entries as well.

We will now define the set $G_n$ by
\[
G_n=\left\{\frac{a}{\left|A_n^{-1}\right|}: a\in\mathbb{Z}
\right\}.
\]
Then, $G_n$ is embedded in $G_{n+1}$ by inclusion, and we write
\[
G=\lim_{\longrightarrow} G_n=\bigcup_{n\ge 1}G_n\subset
\mathbb{Q}.
\]
It is then clear that for
$(\alpha_1,\ldots,\alpha_{n+1})^T\in\mathbb{Z}^{n+1}$,
\[
A_n(\alpha_1,\ldots,\alpha_{n+1})^T\in G_{n+1}.
\]
Giving $G$ the discrete topology, we therefore see that $\Phi_n$
will be a map from $\mathbb{Z}^{n+1}$ to $C(X_{min},G)$.

To show that $K_0(\mathfrak{A})$ is (at least) a subgroup of
$C(X_{min},G)$, we would like to show that the diagram
\[
\xymatrix{ K_0(\mathfrak{A}_n)\cong\mathbb{Z}^{n+1}
\ar[rr]^{\phi_{n*}}
\ar[drr]_{\Phi_n} & & K_0(\mathfrak{A}_{n+1})\cong\mathbb{Z}^{n+2}
 \ar[d]^{\Phi_{n+1}}\\
& & C(X_{min},G)}
\]
commutes, where $\phi_{n\ast}$ is the homomorphism induced by the
standard embedding $\phi_n:\mathfrak{A}_n\to\mathfrak{A}_{n+1}$.

Let $(\alpha_1,\ldots,\alpha_{n+1})^T\in\mathbb{Z}^{n+1}$.  Then,
$\Phi_{n+1}\circ\phi_{n\ast}(\alpha_1,\ldots, \alpha_{n+1})^T$ can
be written as
\[
\left(R_{n+1}\circ [A_{0,1}^{-1}\oplus I_n]\cdots [A_{n-1,n}^{-1}
\oplus I_1]A_{n,n+1}^{-1}\right) \circ
A_{n,n+1}(\alpha_1,\ldots,\alpha_{n+1},0)^T.
\]
If we define $(\beta_1^{n-1},\ldots,\beta_{n+1}^{n-1})^T=
A_{n-1,n}^{-1}(\alpha_1,\ldots,\alpha_{n+1})^T$ and, in general,
for $2\le i<n$,
\[
(\beta_1^{n-i},\ldots,\beta_{n-i+2}^{n-i})^T= A_{n-i,n-i+1}^{-1}
(\beta_1^{n-i+1},\ldots,\beta_{n-i+2}^{n-i+1})^T,
\]
we see that
\begin{eqnarray*}
\lefteqn{\Phi_{n+1}\circ\phi_{n\ast}(\alpha_1,
\ldots,\alpha_{n+1})^T}\\
&=& R_{n+1}\circ [A_{0,1}^{-1}\oplus I_n]\cdots [A_{n-2,n-1}^{-1}
\oplus I_2](\beta_1^{n-1},\ldots,\beta_{n+1}^{n-1},0)^T\\
&=& R_{n+1}(\beta_1^0,\beta_2^0,\beta_3^1,\ldots, \beta_n^{n-2},
\beta_{n+1}^{n-1},0)\\
&=& \beta_1^0\chi_{_{B(1,0)}}+ \sum_{l=1}^{n}
\beta_{l+1}^{l-1}\chi_{_{B(r_l,l)}}.
\end{eqnarray*}
However, the following calculation yields
\begin{eqnarray*}
\Phi_n(\alpha_1,\ldots,\alpha_{n+1})^T
&=&R_n(\beta_1^0,\beta_2^0,\beta_3^1,\ldots,
\beta_n^{n-2},\beta_{n+1}^{n-1})\\
&=& \beta_1^0\chi_{_{B(1,0)}}+ \sum_{l=1}^{n}
\beta_{l+1}^{l-1}\chi_{_{B(r_l,l)}}.
\end{eqnarray*}
Thus, $\Phi_n=\Phi_{n+1}\circ\phi_{n\ast}$, and we conclude that
the diagram commutes.  By the universal property of the direct
limit, it follows that there exists a homomorphism
$\Psi:K_0(\mathfrak{A})\to C(X_{min},G)$.  Since each of the maps
$\Phi_n$ are injective, it follows that in fact, $\Psi$ is an
injection.

\begin{remark}

Considering $C(X_{min},G)$ as an ordered group with positive cone
$C(X_{min},G^+)$, this homomorphism is not necessarily order
preserving. That is, in general,
\[
K_0^+(\mathfrak{A})\not\cong K_0(\mathfrak{A}) \cap
C(X_{min},G^+).
\]
We will say more about this in a moment.  It is true, however,
that the order unit of $K_0(\mathfrak{A})$ in $C(X_{min},G)$ is
$\chi_{_{X_{min}}}$.

\end{remark}

The following theorem provides conditions under which more
information about the structure of $K_0(\mathfrak{A})$ is
available.

\begin{theorem}
\label{rank1}

If, in constructing the matrices $A_{n,n+1}$, integer-valued
columns can be chosen so that $|A_{n,n+1}|=1$, for all $n\ge 0$,
then $K_0(\mathfrak{A})\cong C(X_{min},\mathbb{Z})$ and the order
unit is $\chi_{_{X_{min}}}$.  In general, however,
$K_0^+(\mathfrak{A})\not\cong C(X_{min},\mathbb{Z}^+)$.

\end{theorem}

\begin{proof}

In this case each matrix $A_n$ is invertible and has integer
entries. Thus, $K_0(\mathfrak{A})$ is a subgroup of $C(X_{min},
\mathbb{Z})$. If we let $f\in C(X_{min},\mathbb{Z})$ then by
continuity and the compactness of $X_{min}$, there exists $n\ge 1$
and $(\alpha_1,\ldots,\alpha_{n+1})^T\in\mathbb{Z}^{n+1}$ such
that
\[
f=\sum_{l=1}^{n+1} \alpha_l\chi_{_{B(l,n)}}.
\]
By Lemma \ref{functionform}, there exists
$(\beta_1,\ldots,\beta_{n+1})^T\in\mathbb{Z}^{n+1}$ such that
$f=R_n(\beta_1,\ldots,\beta_{n+1})$.  Thus,
\[
A_n^{-1}(\beta_1,\ldots,\beta_{n+1})^T\in\mathbb{Z}^{n+1} \cong
K_0(\mathfrak{A}_n),
\]
and we have $f=\Phi_n(A_n^{-1}(\beta_1,\ldots,\beta_{n+1})^T)$.

\end{proof}

This theorem implies that whenever two AF algebras satisfy the
hypotheses of Theorem \ref{rank1} and their standard embeddings
can be chosen to yield the same $X_{min}$, then the only aspect of
their dimension groups which distinguishes them is the positive
cone.  This, of course, makes it clear why, in general, $K_0^+
(\mathfrak{A})$ is not $C(X_{min},\mathbb{Z}^+)$.  At the
beginning of this section we commented on how if $X$ is a
$0$-dimensional compact metric space then the AF algebra $C(X)$
has dimension group $(C(X,\mathbb{Z}),C(X,\mathbb{Z}^+),
\chi_{_{X}})$.  Therefore, given $\mathfrak{A}$ satisfying the
hypotheses of Theorem \ref{rank1}, for any $X_{min}$ associated
with $\mathfrak{A}$,
\[
K_0(C(X_{min}))\cong
(C(X_{min},\mathbb{Z}),C(X_{min},\mathbb{Z}^+),
\chi_{_{X_{min}}}).
\]
So, by Elliott's Theorem~\cite{elliott} and Theorem \ref{rank1},
unless $C(X_{min})\cong\mathfrak{A}$, it must be that
$K_0^+(\mathfrak{A})\subset C(X_{min},\mathbb{Z})$ is distinct
from $C(X_{min},\mathbb{Z}^+)$.

One also sees that the flexibility that may exist in choosing
$X_{min}$ is of no help here.  That is, any choice for $X_{min}$
will result in $K_0^+(\mathfrak{A})$ being in general different
from $C(X_{min},\mathbb{Z}^+)$.  Despite
this,~\cite{petersandzerr} shows that useful characterizations of
$K_0^+(\mathfrak{A})$ as a subset of $C(X_{min},\mathbb{Z})$ do
exist in specific cases.  Furthermore, circumstances under which
the topological structure of $X_{min}$ is useful in distinguishing
between AF algebras will be outlined in Theorem \ref{homeomorphic}
below.

Earlier, in showing that $K_0(\mathfrak{A})$ is a subgroup of
$C(X_{min},G)$, we alluded to the fact that it may be strictly
contained in $C(X_{min},G)$.  The following example confirms this.

\begin{example}

Consider the AF algebra $\mathfrak{A}$ with Bratteli diagram
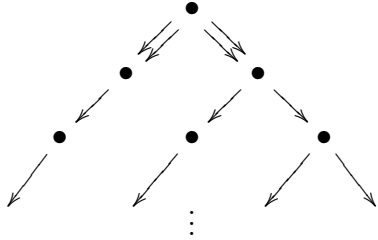
\[
\xymatrix@=12pt{ & & & {\bullet} \ar@<2pt>[dl]\ar@<-2pt>[dl]
\ar@<2pt>[dr]\ar@<-2pt>[dr] & & &\\
  & & {\bullet} \ar[dl] & & {\bullet} \ar[dl] \ar[dr] & &\\
  & {\bullet} \ar[dl]  & & {\bullet} \ar[dl]  & &
  {\bullet} \ar[dl] \ar[dr] &\\
  & & & {\vdots} & & & }
\]
Here, $\displaystyle \overline{A}_{0,1}=\left[\begin{array}{c} 2\\
2\end{array}\right]$ and, for $n\ge 1$,
\[
\overline{A}_{n,n+1}=\left[\begin{array}{ccc} 1 & & \\
 & \ddots & \\
 & & 1\\
 & & 1\end{array}\right],
\]
where all unspecified entries are zero.  Suppose $\displaystyle
A_{0,1}=\left[\begin{array}{cc} 2 & a_{1,2}\\ 2 &
a_{2,2}\end{array} \right]$ and, for $n\ge 1$,
\[
A_{n,n+1}=\left[\begin{array}{cccc} 1 & & & a_{1,n+1}\\
& \ddots & & \vdots\\
& & 1 & a_{n,n+1}\\
& & 1 & a_{n+1,n+1} \end{array}\right],
\]
where the columns are chosen in order to make $A_{n,n+1}$
invertible for all $n\ge 0$. Note that this means
$a_{2,2}-a_{1,2}\ne 0$.  We will show that for all $n\ge 1$, the
vector $\displaystyle \left(
0,\frac{1}{2(a_{2,2}-a_{1,2})},0,\ldots, 0\right)^T\in G_{n+1}$ is
not in the range of $A_n$.  Of course, since $A_n$ is nonsingular,
this is equivalent to showing that
\[
A_n^{-1}\left(0,\frac{1}{2(a_{2,2}-a_{1,2})},0,\ldots,0
\right)^T\notin \mathbb{Z}^{n+1},
\]
for all $n\ge 1$.

It can be easily checked that $\displaystyle
A_n^{-1}\left(0,\frac{1}{2(a_{2,2}-a_{1,2})},0,\ldots,0 \right)^T$
equals
\[
\frac{1}{2(a_{2,2}-a_{1,2})}(a_{1,2},a_{2,2},a_{2,2},\ldots,
a_{2,2})^T.
\]
Now, if $\displaystyle \frac{a_{1,2}}{2(a_{2,2}-a_{1,2})},
\frac{a_{2,2}}{2(a_{2,2}-a_{1,2})}\in\mathbb{Z}$ then
\[
\frac{a_{2,2}}{2(a_{2,2}-a_{1,2})}-\frac{a_{1,2}}{
2(a_{2,2}-a_{1,2})}=\frac{1}{2}
\]
should also be an element of $\mathbb{Z}$.  This contradiction
implies that there exists a function in $C(X_{min},G)$ not in the
range of $\Phi_n$ for any $n$.  Hence, $K_0(\mathfrak{A})$ is a
strict subgroup of $C(X_{min},G)$.

\end{example}

The next important theorem provides conditions under which the
topological structure of $X_{min}$ can be used to help
discriminate between AF algebras.

\begin{theorem}
\label{homeomorphic}

Suppose $\mathfrak{A}_1$ and $\mathfrak{A}_2$ are AF algebras such
that $K_0(\mathfrak{A}_i)\subset C(X_{min}^i,G_i)$ and
$K_0^+(\mathfrak{A}_i)\subset C(X_{min}^i,G_i^+)$, for $i=1,2$,
where both $G_1$ and $G_2$ are subgroups of $\mathbb{Q}$ as in the
discussion preceding Theorem \textup{\ref{rank1}}.  Assume that
for all $U\subset X_{min}^i$ clopen, $\chi_{_{U}}\in
K_0(\mathfrak{A}_i)$. If $\psi:K_0(\mathfrak{A}_1)\to
K_0(\mathfrak{A}_2)$ is an order preserving isomorphism taking
order unit to order unit, then $X_{min}^1$ and $X_{min}^2$ are
homeomorphic.

\end{theorem}

\begin{proof}

To begin, let $x\in X_{min}^1$ and take $\{U_n\}_{n=1}^\infty$ to
be a decreasing sequence of clopen subsets in $X_{min}^1$ with
$\{x\}=\bigcap_{n=1}^\infty U_n$.  Then,
\[
\chi_{_{X_{min}^2}}=\psi(\chi_{_{X_{min}^1}})=\psi(\chi_{_{U_n}})+
\psi(\chi_{_{U_n^c}}),
\]
for all $n\ge 1$, where $U_n^c$ is the complement of $U_n$ in
$X_{min}^1$.  If we suppose the supports of the functions
$\psi(\chi_{_{U_n}})$ and $\psi(\chi_{_{U_n^c}})$ are not
disjoint, then there exists a nonempty clopen set $V\subset
X_{min}^2$ such that $\alpha\chi_{_{V}}\le
\psi(\chi_{_{U_n}}),\psi( \chi_{_{U_n^c}})$, where $\alpha\in
G_2^+$, $\alpha\ne 0$. Because $\psi^{-1}$ is order preserving, it
follows that
\[
\psi^{-1}(\alpha\chi_{_{V}})\le \chi_{_{U_n}},\chi_{_{U_n^c}},
\]
implying $\psi^{-1}(\alpha\chi_{_{V}})=0$ since $U_n\cap U_n^c=
\emptyset$.  Thus, $\alpha\chi_{_{V}}=0$, a contradiction, and we
conclude that $\psi(\chi_{_{U_n}})$ and $\psi(\chi_{_{U_n^c}})$
have disjoint supports.  Therefore, by continuity there exists a
nonempty clopen subset $V_n\subset X_{min}^2$ such that
$\psi(\chi_{_{U_n}})=\chi_{_{V_n}}$, for all $n\ge 1$. Because
$\{\chi_{_{U_n}}\}_{n=1}^\infty$ is a decreasing sequence in
$C(X_{min}^1,G_1)$, $\{\chi_{_{V_n}}\}_{n=1}^\infty$ decreases in
$C(X_{min}^2,G_2)$, and consequently, $\{V_n\}_{n=1}^\infty$ is a
decreasing sequence of clopen sets in $X_{min}^2$.

Let $y\in\bigcap_{n=1}^\infty V_n$ and let $\{W_n\}_{n=1}^\infty$
be a decreasing sequence of clopen sets in $X_{min}^2$ such that
$\{y\}=\bigcap_{n=1}^\infty W_n$ and $W_n\subset V_n$, for all
$n\ge 1$.  Then, $\psi^{-1}(\chi_{_{W_n}})\le \chi_{_{U_n}}$, and
letting $\psi^{-1}(\chi_{_{W_n}})=\chi_{_{Z_n}}$, $Z_n$ clopen, we
have $Z_n\subset U_n$, for all $n\ge 1$.  Furthermore, $x\in Z_n$,
for all $n$ since otherwise $Z_n$ and $U_n$ would eventually be
disjoint.  Hence, $\{x\}=\bigcap_{n=1}^\infty Z_n$.

Now, if $y_1\ne y_2$, $y_1,y_2\in\bigcap_{n=1}^\infty V_n$, then
there would exist sequences $\{W_n^1\}$, $\{W_n^2\}$ of clopen
sets, such that $y_1\in\bigcap_{n=1}^\infty W_n^1$,
$y_2\in\bigcap_{n=1}^\infty W_n^2$, and $W_n^1\cap
W_n^2=\emptyset$, for all $n\ge 1$.  But then
$\psi^{-1}(\chi_{_{W_n^1}})$ and $\psi^{-1}(\chi_{_{W_n^2}})$
would have disjoint supports, and so the corresponding sets
$Z_n^1$, $Z_n^2$, would be disjoint.  This contradicts
$\{x\}=\bigcap_{n=1}^\infty Z_n^i$, $i=1,2$.  Hence,
$\{y\}=\bigcap_{n=1}^\infty V_n$.

Define $f:X_{min}^1\to X_{min}^2$ by letting $f(x)=y$, for all
$x\in X_{min}^1$.  Clearly this map is well defined since any two
sequences $\{U_n\}$, $\{U'_n\}$ which decrease down to $\{x\}$ can
be intertwined, resulting in the intertwining of the corresponding
sequences $\{V_n\}$ and $\{V'_n\}$.  It is also clear that $f$ is
a bijection.  Finally, for $U$ clopen in $X_{min}^1$, one sees
that if $\psi(\chi_{_{U}})=\chi_{_{V}}$ then $f(U)=V$.  Therefore,
$f$ is an open map, and similarly so is its inverse.  Hence,
$X_{min}^1$ is homeomorphic to $X_{min}^2$.

\end{proof}

\begin{remark}

This theorem generalizes the well-known fact that if $C(X_1)$ and
$C(X_2)$ are isomorphic AF algebras, then their spectra $X_1$ and
$X_2$ are homeomorphic.

\end{remark}

To this point, we have only been concerned with those AF algebras
$\displaystyle \mathfrak{A}=\lim_{\longrightarrow}
(\mathfrak{A}_n, \phi_n)$ such that $m_n=n+1$, for all $n\ge 0$.
The case where $m_n=L\in\mathbb{Z}^+$, for all $n\ge 0$ and $L$
fixed, is simpler and more transparent.  For completeness,
however, we mention briefly the key points.

Here, under the assumption that $\overline{A}_{n,n+1}\in
M_L(\mathbb{N})$ is invertible, any minimal Bratteli diagram
corresponding to the set $X_{min}$ will be homeomorphic to a
finite set $\{x_1,\ldots,x_L\}$ with the discrete topology by
Theorem \ref{finitegraphs}.  Then, define the maps
\[
\overline{\Phi}_n=\overline{R}\circ\left[
\overline{A}_{0,1}^{-1}\cdots \overline{A}_{n-1,n}^{-1}\right]:
K_0(\mathfrak{A}_n)\to C(X_{min},G)
\]
where $\overline{R}(\alpha_1,\ldots,\alpha_L)^T=\sum_{i=1}^L
\alpha_i\chi_{_{\{x_i\}}}$,
\[
G_n=\left\{\frac{a}{\left|\overline{A}_{0,1}\right|\cdots
\left|\overline{A}_{n-1,n}\right|}:a\in\mathbb{Z}\right\},
\]
and $\displaystyle G=\lim_{\longrightarrow}G_n$.  It follows that
the diagram
\[
\xymatrix{ K_0(\mathfrak{A}_n)\cong\mathbb{Z}^{L}
\ar[rr]^{\phi_{n*}} \ar[drr]_{\overline{\Phi}_n} & &
K_0(\mathfrak{A}_{n+1})\cong\mathbb{Z}^{L}
 \ar[d]^{\overline{\Phi}_{n+1}}\\
& & C(X_{min},G)}
\]
commutes, and therefore, $K_0(\mathfrak{A})$ is a subgroup of
$C(X_{min},G)$.

As with the earlier discussion, in general, $K_0(\mathfrak{A})$ is
a strict subgroup of $C(X_{min},G)$, the order unit is
$\chi_{_{X_{min}}}$, and the positive cone is not necessarily
$C(X_{min},G^+)$ (or even a subset of it).  We also note here that
Theorems \ref{rank1} and \ref{homeomorphic} remain true in this
context.

\section{Unique Minimal Bratteli Diagrams}
\label{section5}

Considering those AF algebras $\displaystyle
\mathfrak{A}=\lim_{\longrightarrow} (\mathfrak{A}_n,\phi_n)$ such
that $\{m_n\}_{n=0}^\infty$ is strictly monotonically increasing,
we assume now that $\mathfrak{A}$ has a Bratteli diagram with a
unique subgraph corresponding to the set $X_{min}$.  That is to
say, the Bratteli diagram has a unique reduction to a minimal
diagram.  As above, for convenience we dilate the sequence and
assume without a loss of generality that $m_n=n+1$ for all $n\ge
0$.  It is clear that this does not alter the fact that the
diagram corresponding to $X_{min}$ is unique.  One can then show,
up to unitary equivalence, that the connecting homomorphisms
correspond to multiplicity matrices of the form
\[
\overline{A}_{n,n+1}=\left[\begin{array}{cccccc} a_{1,1}^n & & & &
& \\ & \ddots & & & & \\ & & a_{j(n),j(n)}^n & & & \\
& & a_{j(n)+1,j(n)}^n & & & \\ & & & a_{j(n)+2,j(n)+1}^n & & \\
& & & & \ddots & \\ & & & & & a_{n+2,n+1}^n \end{array}\right],
\]
for some $1\le j(n)\le n+1$ and where all unspecified entries are
zero. Therefore, by choosing nonzero integers $b_n$
($=a_{j(n)+1,n+2}$), $n\ge 0$, we can define, for all $n\ge 0$,
\[
A_{n,n+1}=\left[\begin{array}{ccccccc} a_{1,1}^n & & & & &
& \\ & \ddots & & & & & \\ & & a_{j(n),j(n)}^n & & & & \\
& & a_{j(n)+1,j(n)}^n & & & & b_n \\ & & & a_{j(n)+2,j(n)+1}^n & & & \\
& & & & \ddots & & \\ & & & & & a_{n+2,n+1}^n &
\end{array}\right],
\]
in which case $A_{n,n+1}^{-1}$ is the matrix
\[
\left[\begin{array}{ccccccc} {a_{1,1}^n}^{-1} & & & & & &
\\ & \ddots & & & & & \\ & & {a_{j(n),j(n)}^n}^{-1} & 0 & &
& \\ & & & 0 & {a_{j(n)+2,j(n)+1}^n}^{-1} & & \\
& & & & & \ddots & \\ & & & & & & {a_{n+2,n+1}^n}^{-1} \\
& & -\frac{a_{j(n)+1,j(n)}^n}{b_n a_{j(n),j(n)}^n} & \frac{1}{b_n}
& & & 0 \end{array}\right].
\]

Now, letting $(\alpha_1,\alpha_2)^T\in\mathbb{Z}^2$, we see that
\begin{eqnarray*}
\Phi_1(\alpha_1,\alpha_2)^T&=& R_1\left(
\frac{\alpha_1}{a_{1,1}^0}, \frac{\alpha_2}{b_0} -\frac{\alpha_1
a_{2,1}^0}{a_{1,1}^0 b_0}\right)^T\\
&=& \frac{\alpha_1}{a_{1,1}^0}\chi_{_{B(1,0)}}
+\left(\frac{\alpha_2}{b_0} -\frac{\alpha_1 a_{2,1}^0}{a_{1,1}^0
b_0}\right)\chi_{_{B(2,1)}}\\
&=& \frac{\alpha_1}{a_{1,1}^0}\chi_{_{B(1,1)}} +\left(
\frac{\alpha_1}{a_{1,1}^0}+\frac{\alpha_2}{b_0} -\frac{\alpha_1
a_{2,1}^0}{a_{1,1}^0 b_0}\right) \chi_{_{B(2,1)}},
\end{eqnarray*}
since $B(1,0)=B(1,1)\cup B(2,1)$.  Thus,
\[
\Phi_1(\alpha_1,\alpha_2)^T =\frac{\alpha_1}{a_{1,1}^0}
\chi_{_{B(1,1)}} +\frac{\alpha_2}{b_0} \chi_{_{B(2,1)}}
\]
if we choose $b_0=a_{2,1}^0$.  We would like to show that there
exists a sequence $\{b_n\}_{n=0}^\infty$ of positive integers such
that for every $n\ge 1$,
\[
\Phi_n(\alpha_1,\ldots,\alpha_{n+1})^T =\sum_{i=1}^{n+1}
\frac{\alpha_i}{k(i,n)} \chi_{_{B(i,n)}},
\]
where $k(i,n)\in\mathbb{Z}^+$, for all $1\le i\le n+1$.  Our above
calculations prove that it is possible to choose $b_0$
accordingly.

Assume that $b_0,\ldots, b_{k-1}\in\mathbb{Z}^+$ have been chosen
in order to guarantee $\Phi_k$ is of the appropriate form for
$1\le k\le n$.  Let $(\alpha_1,\ldots, \alpha_{n+2})^T\in
\mathbb{Z}^{n+2}$.  Then,
\[
\Phi_{n+1}(\alpha_1,\ldots, \alpha_{n+2})^T =R_{n+1}
\left[\begin{array}{cc} A_n & 0\\ 0 & 1\end{array}\right]
A_{n,n+1}^{-1}(\alpha_1,\ldots,\alpha_{n+2})^T.
\]
By then first computing
$A_{n,n+1}^{-1}(\alpha_1,\ldots,\alpha_{n+2})^T$, one is able to
show that $\Phi_{n+1}(\alpha_1,\ldots,\alpha_{n+2})^T$ is equal to
\begin{eqnarray*}
\mbox{} &\mbox{}& \Phi_n \left( \frac{\alpha_1}{a_{1,1}^n},\ldots,
\frac{\alpha_{j(n)}}{ a_{j(n),j(n)}^n},
\frac{\alpha_{j(n)+2}}{a_{j(n)+2,j(n)+1}^n},
\ldots, \frac{\alpha_{n+2}}{a_{n+2,n+1}^n}\right)\\
&+& \left(\frac{\alpha_{j(n)+1}}{b_n}
-\frac{\alpha_{j(n)}a_{j(n)+1,j(n)}^n}{ b_n a_{j(n),j(n)}^n}
\right)\chi_{_{B(r_{n+1},n+1)}}.
\end{eqnarray*}
Hence, our induction hypothesis allows us to write
$\Phi_{n+1}(\alpha_1,\ldots, \alpha_{n+2})^T$ as
\begin{eqnarray*}
\mbox{}&\mbox{}& \sum_{i=1}^{j(n)} \frac{\alpha_i}{a_{i,i}^n
k(i,n)} \chi_{_{B(i,n)}} +\sum_{i=j(n)+1}^{n+1}
\frac{\alpha_{i+1}}{a_{i+1,i}^n k(i,n)}\chi_{_{B(i,n)}} \\
&+&\left(\frac{\alpha_{j(n)+1}}{b_n}
-\frac{\alpha_{j(n)}a_{j(n)+1,j(n)}^n}{ b_n a_{j(n),j(n)}^n}
\right)\chi_{_{B(r_{n+1},n+1)}}.
\end{eqnarray*}
Now, we note that for $1\le i\le j(n)-1$, $B(i,n)=B(i,n+1)$, and
for $j(n)+1\le i\le n+1$, $B(i,n)=B(i+1,n+1)$.  Furthermore,
$B(j(n),n)=B(j(n),n+1)\cup B(j(n)+1,n+1)$ and $r_{n+1}=j(n)+1$.
Thus, $\Phi_{n+1}(\alpha_1,\ldots,\alpha_{n+2})^T$ becomes
\begin{eqnarray*}
\mbox{} &\mbox{}& \sum_{i=1}^{j(n)-1} \frac{\alpha_i}{ a_{i,i}^n
k(i,n)}\chi_{_{B(i,n+1)}} +\sum_{i=j(n)+1}^{n+1}
\frac{\alpha_{i+1}}{a_{i+1,i}^n k(i,n)}\chi_{_{B(i+1,n+1)}}\\
&+&\frac{\alpha_{j(n)}}{a_{j(n),j(n)}^n k(j(n),n)} \left(
\chi_{_{B(j(n),n+1)}}+\chi_{_{B(j(n)+1,n+1)}}\right)\\
&+&\left(\frac{\alpha_{j(n)+1}}{b_n} -\frac{\alpha_{j(n)}
a_{j(n)+1,j(n)}^n}{b_n a_{j(n),j(n)}^n}\right)
\chi_{_{B(j(n)+1,n+1)}}\\
&=&\sum_{i=1}^{j(n)} \frac{\alpha_i}{ a_{i,i}^n
k(i,n)}\chi_{_{B(i,n+1)}} +\sum_{i=j(n)+1}^{n+1}
\frac{\alpha_{i+1}}{a_{i+1,i}^n k(i,n)}\chi_{_{B(i+1,n+1)}}\\
&+& \left[ \frac{\alpha_{j(n)}}{a_{j(n),j(n)}^n k(j(n),n)} \left(
1-\frac{a_{j(n)+1,j(n)}^n k(j(n),n)}{b_n}\right)\right]
\chi_{_{B(j(n)+1,n+1)}}\\
&+&\frac{\alpha_{j(n)+1}}{b_n}\chi_{_{B(j(n)+1,n+1)}}.
\end{eqnarray*}
If we let $b_n=a_{j(n)+1,j(n)}^n k(j(n),n)$, then we have the
desired outcome.

The significance of this is that when determining
$K_0^+(\mathfrak{A})$, we see that since
$K_0^+(\mathfrak{A})=\bigcup_{n\ge 1}\Phi_n(\mathbb{Z}_+^{n+1})$,
clearly $K_0^+(\mathfrak{A})\subset C(X_{min},G^+)$.  Furthermore,
for all $U\subset X_{min}$ clopen, $\chi_{_{U}}\in
K_0(\mathfrak{A})$.  Therefore, the hypotheses of Theorem
\ref{homeomorphic} are satisfied.  We then immediately have the
following corollaries.

\begin{corollary}
\label{cor}

Let $\mathfrak{A}_1$ and $\mathfrak{A}_2$ be AF algebras, each
which have Bratteli diagrams for which there are unique minimal
reductions corresponding to the sets $X_{min}^1$ and $X_{min}^2$.
If $\mathfrak{A}_1 \cong \mathfrak{A}_2$ then $X_{min}^1$ is
homeomorphic to $X_{min}^2$.

\end{corollary}

\begin{corollary}
\label{importantcor}

Suppose $\mathfrak{A}$ is an AF algebra with a Bratteli diagram
for which there is a unique minimal reduction corresponding to the
set $X_{min}$. Then, every order preserving automorphism $\psi$ of
$K_0(\mathfrak{A})$ which takes the order unit onto itself must be
of the form $\psi(f)=f\circ{\theta}^{-1}$ for some homeomorphism
$\theta$ of $X_{min}$.

\end{corollary}

\begin{remark}

If $X$ is a $0$-dimensional compact metric space then the above
arguments show that with $b_n=1$, for all $n\ge 0$, the methods of
this paper provide a direct generalization of the usual way of
showing that
\[
K_0(C(X))\cong\left( C(X,\mathbb{Z}),C(X,\mathbb{Z}^+),
\chi_{_{X}}\right).
\]

\end{remark}

\begin{remark}

When $\mathfrak{A}$ is such that $m_n=L<\infty$, for all $n\ge 1$,
the existence of a unique minimal Bratteli diagram implies the
multiplicity matrices are permutation similar to diagonal
matrices.  Thus, these corollaries generalize to this situation as
well.

\end{remark}

We conclude this paper with some applications of Corollary
\ref{importantcor} which provide information about the
automorphism groups of certain AF algebras.

\begin{corollary}
\label{autcor}

Suppose $\mathfrak{A}$ is an AF algebra with a Bratteli diagram
for which there is a unique minimal reduction corresponding to the
set $X_{min}$. If the only homeomorphism of $X_{min}$ is the
identity map, then $\mbox{Aut}(\mathfrak{A})=
\overline{\mbox{Inn}(\mathfrak{A})}$.

\end{corollary}

\begin{proof}

Each automorphism $\alpha\in\mbox{Aut}(\mathfrak{A})$ induces an
automorphism $\alpha_{\ast}\in\mbox{Aut}(K_0(\mathfrak{A}))$.  By
Corollary \ref{importantcor}, $\alpha_{\ast}(g)=g\circ\theta^{-1}$
for some homeomorphism $\theta$ of $X_{min}$, where we are
representing $K_0(\mathfrak{A})$ as a subset of $C(X_{min},G)$. By
hypothesis, it follows that $\alpha_{\ast}=\mbox{id}_{\ast}$, and
therefore, by~\cite[Theorem IV.5.7]{davidson},
$\alpha\in\overline{\mbox{Inn}(\mathfrak{A})}$. Hence,
$\mbox{Aut}(\mathfrak{A})=\overline{\mbox{Inn}(\mathfrak{A})}$.

\end{proof}

\begin{remark}

In particular, this applies to the UHF algebras, and should be
compared with~\cite[Corollary IV.5.8]{davidson}.

\end{remark}

The converse of Corollary \ref{autcor} is in general not true. The
following example illustrates this and demonstrates an application
of Corollary \ref{importantcor} which utilizes the structure of
$X_{min}$ to obtain complete information about the automorphism
group of $K_0(\mathfrak{A})$.

\begin{example}

Consider the AF algebra $\mathfrak{A}$ with Bratteli diagram
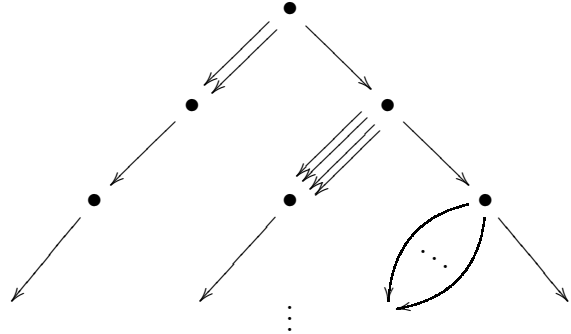
\[
\xymatrix{ & & & {\bullet} \ar@<2pt>[dl]\ar@<-2pt>[dl] \ar[dr] & & &\\
  & & {\bullet} \ar[dl] & & {\bullet} \ar@<5pt>[dl]\ar@<-5pt>[dl]
  \ar@<2pt>[dl]\ar@<-2pt>[dl] \ar[dr] & &\\
  & {\bullet} \ar[dl]  & & {\bullet} \ar[dl]  & &
  {\bullet} \ar@/^1pc/[dl]_{\ddots}\ar@/_1pc/[dl] \ar[dr] &\\
  & & & {\vdots} & & & }
\]
where the multiplicity matrices are such that
\[
\overline{A}_{n,n+1}=\left[\begin{array}{cccc} 1 & &
& \\ & \ddots & & \\ & & 1 & \\
& & & 2^{n+1} \\ & & & 1 \end{array}\right]\in M_{n+2,n+1},
\]
for all $n\ge 0$.  We then complete these matrices in such a way
so as to obtain
\[
A_{n,n+1}=\left[\begin{array}{ccccc} 1 & & & & \\ & \ddots & & &
\\ & & 1 & & \\ & & & 2^{n+1} & \\ & & & 1 & 1\end{array}\right]
\in M_{n+2}
\]
for all $n\ge 0$.  Therefore, if we define $A_n$ as in Section
\ref{sec4}, one can verify by induction on $n$ that
\[
A_n=\left[\begin{array}{ccccc} 2^{-1} & & & & \\ -2^{-1} & 2^{-2}
& & & \\ & -2^{-2} & & & \\ & & \ddots & & \\ & & & 2^{-n} & \\ &
& & -2^{-n} & 1\end{array}\right]
\]
for all $n\ge 1$.

Now, with $X_{min}$ given (uniquely) as the set of all infinite
paths in the graph
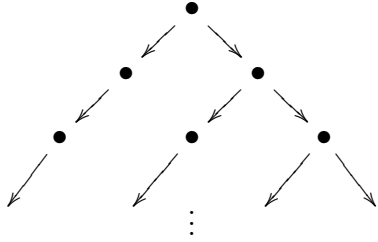
\[
\xymatrix@=12pt{ & & & {\bullet} \ar[dl] \ar[dr] & & &\\
  & & {\bullet} \ar[dl] & & {\bullet} \ar[dl]\ar[dr] & &\\
  & {\bullet} \ar[dl]  & & {\bullet} \ar[dl]  & &
  {\bullet} \ar[dl] \ar[dr] &\\
  & & & {\vdots} & & & }
\]
we have, for any $(\alpha_1,\ldots,\alpha_{n+1})^T\in
\mathbb{Z}^{n+1}$, that $\Phi_n(\alpha_1,\ldots,\alpha_{n+1})^T$
is equal to
\[
\alpha_1 2^{-1}\chi_{_{B(1,0)}} +\sum_{l=1}^{n-1} (\alpha_{l+1}
2^{-(l+1)} -\alpha_l 2^{-l})\chi_{_{B(l+1,l)}}
+(\alpha_{n+1}-\alpha_n 2^{-n})\chi_{_{B(n+1,n)}},
\]
which, due to the structure of $X_{min}$, can be written as
\[
\left(\sum_{l=1}^n \alpha_l 2^{-l}\chi_{_{B(l,n)}}\right)
+\alpha_{n+1} \chi_{_{B(n+1,n)}}.
\]

Now, by Corollary \ref{importantcor}, every order preserving
automorphism $\psi$ of $K_0(\mathfrak{A})$ must be of the form
$\psi(f)=f\circ \theta^{-1}$ for some homeomorphism $\theta$ of
$X_{min}$, where $f\in K_0(\mathfrak{A})\subset C(X_{min},G)$. The
set $X_{min}$ is easily seen to be homeomorphic to
\[
\overline{\left\{\frac{1}{n}:n\in\mathbb{Z}^+\right\}},
\]
and thus, any homeomorphism $\theta$ of $X_{min}$ must fix the
point corresponding to $0$.  It follows that if $\theta$ is a
non-identity homeomorphism of $X_{min}$, then there exists $n\ge
1$ and $1\le l_1<l_2\le n$ such that
\[
\theta(B(l_1,n))=B(l_2,n),
\]
where here, $B(l_1,n)$ and $B(l_2,n)$ are singletons.  But then,
for the function
$f=2^{-l_1}\chi_{_{B(l_1,n)}}+2^{-l_2}\chi_{_{B(l_2,n)}}\in
K_0(\mathfrak{A})$, we see that
\[
\psi(f)=f\circ\theta^{-1} =2^{-l_1}\chi_{_{B(l_2,n)}} +2^{-l_2}
\chi_{_{B(l_1,n)}}\notin K_0(\mathfrak{A}).
\]
Hence, we conclude that the only order preserving automorphism of
$K_0(\mathfrak{A})$ is the identity map.

Consequently, the converse of Corollary \ref{autcor} is not true.
After all, we have just shown that the only order preserving
automorphism of $K_0(\mathfrak{A})$ is the identity map.  So
again, by~\cite[Theorem IV.5.7]{davidson}, $Aut(\mathfrak{A})
=\overline{Inn(\mathfrak{A})}$.  However, $X_{min}$ clearly has
many nontrivial homeomorphisms.

\end{example}

\begin{corollary}
\label{cor2}

If $X$ is a $0$-dimensional compact metric space then
$\mbox{Aut}(C(X))=\{\widehat{\theta}:\theta\mbox{ is a
homeomorphism of }X\}$ where $\widehat{\theta}(f)=f\circ
\theta^{-1}$.

\end{corollary}

\begin{proof}

Clearly, if $\theta$ is a homeomorphism on $X$, then
$\widehat{\theta}$ is an automorphism of $C(X)$.  Thus, we must
show all automorphisms are of this form.

Let $\alpha\in\mbox{Aut}(C(X))$.  Then, $\alpha$ induces an
automorphism $\alpha_{\ast}$ of $K_0(C(X))$, which by Corollary
\ref{importantcor} must be of the form $\alpha_{\ast}(f)=f\circ
\theta^{-1}$ for some homeomorphism $\theta$ of $X$.  But then,
for all $f\in K_0(\mathfrak{A})$,
\[
{\widehat{\theta}}_{\ast}^{-1}\circ\alpha_{\ast}(f)
={\widehat{\theta}}_{\ast}^{-1}\left( f\circ\theta^{-1}\right)
=f\circ\theta^{-1}\circ\theta =f.
\]
Hence, ${\widehat{\theta}}_{\ast}^{-1}\circ\alpha_{\ast}
=\mbox{id}_{\ast}$.

By~\cite[Theorem IV.5.7]{davidson}, it follows that
${\widehat{\theta}}^{-1}\circ \alpha$ is an approximately inner
automorphism of $C(X)$, which of course means
${\widehat{\theta}}^{-1}\circ\alpha =\mbox{id}$. Therefore,
$\alpha=\widehat{\theta}$.

\end{proof}

Possibly a more interesting observation than these last two
results is that Corollary \ref{importantcor} in some sense
generalizes Corollary \ref{cor2} to a larger class of algebras.
One surprising aspect of this is that the correct set of
homeomorphisms is not on the spectrum of $\mathfrak{D}\subset
\mathfrak{A}$, but rather on the (often) smaller set $X_{min}$.

Since the results of this section are for AF algebras with
Bratteli diagrams that have unique reductions to $X_{min}$, we
state the following result which provides a characterization of
these algebras.

\begin{theorem}

Suppose there exists a way to complete the matrices
$\{\overline{A}_{n,n+1}\}_{n=0}^\infty$ such that for all $n\ge
1$,
\[
\Phi_n(\alpha_1,\ldots,\alpha_{n+1})^T =\sum_{i=1}^{n+1}
\frac{\alpha_i}{k(i,n)}\chi_{_{B(i,n)}},
\]
where $k(i,n)\in\mathbb{Z}^+$ for all $1\le i\le n+1$.  Then, any
two choices for the set $X_{min}$ must be homeomorphic.

\end{theorem}

\begin{proof}

In this case, for any $X_{min}$ associated to $\mathfrak{A}$, we
have $K_0^+(\mathfrak{A})=K_0(\mathfrak{A})\cap C(X_{min},G^+)$
and for all $U\subset X_{min}$, clopen, the function $\chi_{_{U}}$
is an element of $K_0(\mathfrak{A})$.  Suppose $X_{min}^1$ and
$X_{min}^2$ are two choices for the set $X_{min}$.  We then know
by Elliott's Theorem~\cite{elliott} that there exists an order
preserving isomorphism
\[
\psi:K_0(\mathfrak{A})\subset C(X_{min}^1,G_1)\to
K_0(\mathfrak{A})\subset C(X_{min}^2,G_2)
\]
which takes order unit to order unit.  By Theorem
\ref{homeomorphic}, it follows that $X_{min}^1$ and $X_{min}^2$
are homeomorphic.

\end{proof}

\bibliography{author}

\end{document}